\input amstex

\input epsf
\epsfverbosetrue \input amsppt.sty
\magnification 980\vsize=21 true cm \hsize=16 true cm \voffset=1.1
true cm \pageno=1 \NoRunningHeads \TagsOnRight
\def\p{\partial}
\def\ve{\varepsilon}
\def\f{\frac}

\def\na{\nabla}
\def\la{\lambda}
\def\al{\alpha}

\def\g{\gamma}
\def\G{\Gamma}

\def\dl{\delta}

\def\q{\qquad}

\def\ds{\displaystyle}


\topmatter


\title Blowup of classical solutions for a class of 3-D quasilinear wave
equations with small initial data
\endtitle

\author Bingbing Ding, Ingo Witt, and Huicheng Yin\endauthor

\address Department of Mathematics and IMS, Nanjing University,
Nanjing 210093, P.R.~China \endaddress


\email balancenjust\@yahoo.com.cn \endemail

\thanks Ding Bingbing and Yin Huicheng were
supported by the NSFC (No.~10931007, No.~11025105) and by the
Priority Academic Program Development of Jiangsu Higher Education
Institutions. Ingo Witt was supported by the DFG via the
Sino-German project ``Analysis of PDEs and Applications.''
This research was carried out when Ingo Witt was visiting Nanjing
University in March of 2012.
\endthanks

\address Mathematical Institute, University of G\"ottingen,
Bunsenstr.~3-5, D-37073 G\"ottingen, Germany \endaddress

\email iwitt\@uni-math.gwdg.de \endemail

\address Department of Mathematics and IMS, Nanjing University,
Nanjing 210093, P.R.~China \endaddress

\email huicheng\@nju.edu.cn \endemail




\abstract
This paper is concerned with the small smooth data problem for the 3-D
nonlinear wave equation $\p_t^2u-\left(1+u+\p_t u\right)\Delta
u=0$. This equation is prototypical of the more general equation
$\dsize\sum_{i,j=0}^3g_{ij}(u, \na u)\p_{ij}u=0$, where $x_0=t$ and
$g_{ij}(u, \na
u)=c_{ij}+d_{ij}u+\ds\sum_{k=0}^3e_{ij}^k\p_ku+O(|u|^2+|\na u|^2)$ are
smooth functions of their arguments, with $c_{ij}, d_{ij}$ and
$e_{ij}^k$ being constants, and $d_{ij}\neq0$ for some $(i,j)$;
moreover, $\dsize\sum_{i,j,k=0}^3e_{ij}^k(\p_ku)\p_{ij} u$ does not
fulfil the null condition. For the 3-D nonlinear wave equations
$\p_t^2u-\left(1+u\right)\Delta u=0$ and $\p_t^2u-\left(1+\p_t
u\right)\Delta u=0$, H.~Lindblad, S.~Alinhac, and F.~John proved and
disproved, respectively, the global existence of small smooth data
solutions. For radial initial data, we show that the small smooth data
solution of $\p_t^2u-\left(1+u+\p_t u\right)\Delta u=0$ blows up in
finite time. The explicit expression of the asymptotic lifespan
$T_{\ve}$ as $\ve\to0^+$ is also given.
\endabstract


\endtopmatter


\document

\head \S 1. Introduction and main results \endhead

\vskip 0.3 true cm

We consider the second-order nonlinear wave equation in
$[0, \infty)\times \Bbb R^n$
$$
\cases
&\dsize\sum_{i,j=0}^ng_{ij}(u, \na
u)\p_{ij}u=0,\\
&\left(u(0,x), \p_t u(0,x)\right)=\left(\ve u_0(x), \ve u_1(x)\right),
\endcases\tag1.1
$$
where $x_0=t$, $x=(x_1, ..., x_n)$, $\na=(\p_0, \p_1, ... , \p_n)$,
$\ve>0$ is a sufficiently small constant, $u_0(x), u_1(x)\in
C_0^{\infty}(\Bbb R^n)$, and the $g_{ij}(u, \na u)$ are smooth functions of
there arguments which are of the form
$$g_{ij}(u, \na
u)=c_{ij}+d_{ij}u+\ds\sum_{k=0}^ne_{ij}^k\p_ku+O(|u|^2+|\na
u|^2)\tag1.2$$ with $c_{ij}, d_{ij}$ and $e_{ij}^k$ being
constants. We assume that the linear part
$\dsize\sum_{i,j=0}^nc_{ij}\p_{ij}u$ is strictly hyperbolic with
respect to time $t$. From [8-9, 14-16] we have that (1.1)
has a global smooth solution when $n\ge 4$.

If $n=3$ and $d_{ij}=0$ for all $0\le i, j\le 3$ in (1.2), then
(1.1) has a global smooth solution if the null condition for the the
main part $\dsize\sum_{i,j,k=0}^3e_{ij}^k\p_ku\p_{ij}u$ holds,
otherwise the solution of (1.1) blows up in finite time.  See [4,
8-13, 17, 20-27] and the references therein.

If $n=3$ and $d_{ij}\not=0$ for some $(i, j)$, but $e_{ij}^k=0$ for all
$0\le i, j, k\le 3$ in (1.2), then it follows from the results in [3]
and [18-19] that (1.1) has a global smooth solution.

The following interesting problem naturally arises: If $n=3$,
$d_{ij}\not=0$ for some $(i, j)$, and
$\dsize\sum_{i,j,k=0}^3e_{ij}^k\p_ku\p_{ij}u$ in (1.2) does not
fulfill the null condition, does the smooth solution of (1.1) blow up
in finite time or does it exist globally? In this paper, we are
concerned with this problem, especially (and without loss of
generality) the protypical equation $\p_t^2u-(1+u+\p_t u)\Delta u=0$
is studied. More specifically, we consider the problem
$$\cases
&\p_t^2u-(1+u+\p_tu)\triangle u=0,\quad\text{$(t,x)\in [0,\infty)\times\Bbb R^3$},\\
&(u(0,x), \p_tu(0,x))=(\ve u_0(x),\ve u_1(x)),\\
\endcases\tag1.3
$$
where $u_0(x), u_1(x)\in C_0^\infty(\Bbb R^3)$ are radially
symmetric and
$\operatorname{supp}u_0\cup \operatorname{supp}u_1\subseteq\{x :
|x|\leq M\}$ with $M>0$ a constant. For notational convenience, we
write $(u_0(r), u_1(r))$ instead of $(u_0(x), u_1(x))$ later on and
the domains of definition of $u_0(r)$ and $u_1(r)$ are simultaneously
extended to $[-M, M]$. This results from the fact that $u_0(r)$ and
$u_1(r)$ are actually smooth functions of $r^2$ due to $u_0(x),
u_1(x)\in C_0^\infty(\Bbb R^3)$.

Let $F_0(s)=\ds\f12\biggl(su_0(s)+\int_s^{\infty} su_1(s)\,ds\biggr)$
for $s\in \Bbb R$.  According to Theorem 6.2.2 and (6.2.12) of [9],
we know that the function $F_0(s)\not\equiv 0$ unless both
$u_0(s)\equiv 0$ and $u_1(s)\equiv 0$. Moreover, $F_0(s)\equiv 0$ for
$|s|\ge M$.

Let $\tau(s)=\ds \f{2}{F_0'(s)}\ln\f{F_0''(s)}{F_0''(s)-F_0'(s)}$ for
$s\in [-M, M]$, $F_0'(s)\not =0$, and
$\ds\f{F_0''(s)}{F_0''(s)-F_0'(s)}>0$. Further let \linebreak $A=\biggl\{s\in
(-M, M): F_0'(s)\not
=0, \ds\f{F_0''(s)}{F_0''(s)-F_0'(s)}>0, \tau(s)>0\biggr\}$ and
$B=\{s\in (-M, M): F_0'(s)=0, F_0''(s)>0\}$, and denote
$$
\tau_0=\min\biggl\{\ds\min_{s\in A}\tau(s),
\ds\min_{s\in B}\ds\f{2}{F_0''(s)}\biggr\}.\tag1.4
$$
It can be shown that $\tau_0$ is a finite positive number
if $(u_0(r), u_1(r))\not\equiv 0$ holds.

The main result of this paper is:

\proclaim{Theorem 1.1} Assume that $u_0(x)$, $u_1(x)\in
C_0^{\infty}(\Bbb R^3)$ only depend on
$r=\sqrt{x_1^2+x_2^2+x_3^2}$. If $u_0(x)\not\equiv 0$ or
$u_1(x)\not\equiv 0$, then\/ \rom{problem (1.3)} has a
$C^\infty$ solution $u(t,x)$ for $0\leq t <T_\varepsilon$, where $T_{\ve}$
stands for the lifespan of the smooth solution $u(t,x)$, which
satisfies
$$
\lim_{\varepsilon\rightarrow 0}\varepsilon\ln T_\varepsilon =\tau _0>0.\tag1.5
$$
\endproclaim

\remark{Remark 1.1}
It follows from Theorem 1.1 that the smooth solution of (1.3) blows up
in finite time provided that $(u_0(r), u_1(r))\not\equiv 0$.
\endremark

\remark{Remark 1.2}
We have asserted that $\tau_0>0$ is a finite number as long as
$(u_0(r), u_1(r))\not\equiv 0$. To prove this, it suffices to show
that $\tau_1=\ds\min_{s\in A}\tau(s)>0$ is finite, since
$\tau_2=\ds\min_{s\in B}\ds\f{2}{F_0''(s)}>0$ obviously holds. We
first show $A\not=\emptyset$. Let $A_1=\left\{s\in (-M, M):
F_0''(s)>0\right\}$.  Then obviously $A\subseteq A_1$. Since $A_1$ is
a bounded, open, and nonempty set by $(u_0(r), u_1(r))\not\equiv 0$
and the continuity of $F_0''(s)$, one can write $A_1=\ds\bigcup_l(a_l,
b_l)$, where the union is finite or countable infinite and $(a_{l_1},
b_{l_1})$ and $(a_{l_2}, b_{l_2})$ are disjoint for $l_1\not=l_2$. Moreover,
$F''_0(a_l)=F''_0(b_l)=0$ for any $l$. Set $B_0=\left\{s\in [-M,
M]: F_0'(s)=0\right\}$. Then $B_0$ is closed and $(a_l, b_l)\setminus B_0$
is an open and nonempty set.  Write $(a_l, b_l)\setminus B_0=
\ds\bigcup_m(a_{lm}, b_{lm})$, where $F_0'(a_{lm})=F_0'(b_{lm})=0$ and the intervals
$(a_{lm}, b_{lm})$ are disjoint for different $(l,m)$. For any $(l_0,
m_0)$, there exists at least exists one point $s_0\in \left(a_{l_0m_0},
b_{l_0m_0}\right)$  such that
$F_0''(s_0)-F_0'(s_0)>0$. If not, then $F_0''(s)-F_0'(s)\le 0$ for
$s\in (a_{l_0m_0}, b_{l_0m_0})$, which means $(e^{-s}F_0'(s))'\le 0$
and then $F_0'(s)\equiv 0$ because of $F_0'(a_{lm})=F_0'(b_{lm})=0$.
This contradicts $F_0''(s)>0$ for $s\in
\left(a_{l_0m_0}, b_{l_0m_0}\right)\subseteq A_1$. On the other hand, $\tau (s_0)>0$
obviously holds, and thus the set $A$ is nonempty. Next we show
$\tau_1>0$. Set $z=-\ds\f{F_0'(s)}{F_0''(s)}$ for $s\in A$. Then
$\tau(s)=\ds\f{2}{F_0''(s)}\ds\f{\ln(1+z)}{z}$ with $s\in A$, $z>-1$
and $z\not=0$.  If $z\in (-1, 0)$, then
$\tau(s)\ge\ds\f{1}{\ds\max_{s\in A} F_0''(s)}$. If $z\in (0, N]$ for
some fixed $N>0$, then $\tau(s)\ge\ds\f{C_N}{\ds\max_{s\in A}
F_0''(s)}$. If $z>0$ is large, then $F_0''(s)>0$ is small and
$F_0'(s)<0$ holds. In this case, we conclude from
$\tau(s)=-\ds\f{2}{F_0'(s)}\ln(1+z)$ that, if $|F_0'(s)|\ge \dl_0$ for
some fixed constant $\dl_0>0$, then $\tau(s)\ge C(\dl_0)>0$; on the
other hand, if $|F_0'(s)|$ is small, then $\tau(s)>0$ is
large. Consequently, putting everything together, the constants
$\tau_1$ and $\tau_0$ are positive and finite.
\endremark

\remark{Remark 1.3} The result of this paper can be extended to the more
general nonlinear wave equation
$$
\cases
&\dsize\sum_{i,j=0}^3g_{ij}(u, \na u)\p_{ij}u=0,\\ &(u(0,x), \p_t
u(0,x))=(\ve u_0(x), \ve u_1(x)),
\endcases\tag1.6
$$
where $x_0=t$, $u_0(x), u_1(x)\in C_0^{\infty}(\Bbb R^n)$,
the $g_{ij}(u, \na u)$ are smooth functions of their arguments which are of the form
$$g_{ij}(u, \na
u)=c_{ij}+d_{ij}u+\ds\sum_{k=0}^3e_{ij}^k\p_ku+O(|u|^2+|\na u|^2)$$
with $c_{ij}, d_{ij}$ and $e_{ij}^k$ being constants, and
$d_{ij}\not=0$ for some $(i,j)$, and the null condition for
$\dsize\sum_{i,j,k=0}^3e_{ij}^k\p_ku\p_{ij}^2u$ does not hold.
Because a proof of this statement just requires the methods and
estimates used in this paper and the ones of [5] and [25] combined
with blowup system techniques of [1-2], but it is technical and
tendious otherwise, it is omitted here.
\endremark

\remark{Remark 1.4} For the 2-D case of problem (1.1),
if the coefficients $g_{ij}(u, \na u)$ are independent of $u$, then
there is a rather complete collection of results on the global
existence and the blowup, respectively, of small smooth data
solutions, see [1-2, 7] and the references therein. On the contrary,
if the coefficients $g_{ij}(u,\na u)$ depend on both $u$ and $\na u$,
there have been no systematic studies so far. Related results will
appear in a forthcoming paper of ours.
\endremark

\smallskip

In order to prove Theorem 1.1, we first derive the lower bound on the
lifespan $T_{\ve}$ for problem (1.3) when the initial data are radial.
By constructing an approximate solution as in [9] or [6], then
considering the difference of the exact solution and the approximate
solution, applying the Klainerman-Sobolev inequality, and establishing
some further energy estimates, we obtain a lower bound on the
lifespan $T_{\ve}$. Here we point out that although $\ve \ln
T_{\ve}\ge C>0$ has already been shown in [17], for the reader's
convenience and to obtain the sharp lower bound $\tau_0$, we still
give a complete proof.  On the other hand, it follows from radial
symmetry of the initial data $(u_0(x), u_1(x))$ that the solution
$u(t,x)$ is also radially symmetric for $t<T_{\ve}$. Based on this, we
change (1.3) into a $2\times 2$ system of two independent variables
$(t, r)$. Then, by using the properties of the approximate solution
constructed above and some delicate analysis, and by treating the
solution $u$ accordingly, we obtain the upper bound on $T_{\ve}$. Here
the derivation is motivated by the methods of [10], where the equation
$\p_t^2u-c^2(\p_tu)\Delta u=0$ was studied, the coefficient $c^2(\p_tu)$
of which only depends on the gradient of the solution $u$. Compared
with [10] and [6], due to the simultaneous appearance of $u$ and $\p_t
u$ in the coefficients of equation (1.3), we have to introduce a few
more quantities in order to get a ``blowup''-type nonlinear
second-order ordinary differential equation with suitable initial data
that provides the upper bound on $\ve \ln T_{\ve}$.  Based on the
results in the two steps above, we finally complete the proof of
Theorem 1.1.

\smallskip

In this paper, we will use the following notation:

$Z$ denotes one of the Klainerman vector fields in the radially
symmetric case, i.e.,
$$
\p_r, \ \p_t, \ S=t\p_t+r\p_r, \ H=r\p_t+t\p_r,
$$
$\p$ stand for $\p_r$ or $\p_t$, and the norm $\|f\|_{L^2}$ means
$\|f(t,\cdot)\|_{L^2(\Bbb R^3)}$.


\head \S 2. The lower bound on the lifespan $T_\varepsilon$ \endhead

In this section, we establish the lower bound on $T_\varepsilon$ for
smooth solutions of the Cauchy problem (1.3). As in the proof of [9,
Theorem 6.5.3], by constructing an approximate solution $u_a$ of (1.3)
and then estimating the difference of $u_{a}$ and the solution $u$, we
obtain the lower bound on $T_{\ve}$ by a continuity induction
argument. The new ingredients in this procedure are the construction
of the approximate solution and treating the solution $u$ itself that
occurs in the equation in (1.3) rather than the derivatives of this
solution only as in [9].  Although this procedure is analogous to the
one in [6], for the reader's convenience and also as it is applied to
obtain the upper bound on $T_{\ve}$, we still give a complete proof.

Let the slow time variable be $\tau=\varepsilon \ln(1+t)$, and assume the
solution of (1.3) is approximated by
$$
\f{\varepsilon}{r}\,V(q,\tau),\quad r>0,
$$
where $q=r-t$ and $V(q, \tau)$ solves the equation
$$\cases
2\partial_{q\tau} V+V\p_q^2 V-\p_q V\p_q^2V=0,
\qquad\text{$(q,\tau)\in \Bbb R\times[0,\infty)$} ,\\
V(q,0)=F_0(q) ,\\
\operatorname{supp}V\subseteq \{q\leq M\}.\\
\endcases\tag2.1
$$


\proclaim{Lemma 2.1} \rom{(2.1)} has a $C^\infty$ solution for $0\leq\tau<\tau_0$,
where $\tau_0$ is given by \rom{(1.4)}.
\endproclaim

\demo{Proof} Set $w(q,\tau)=\partial_q V(q,\tau)$. Then it follows from (2.1) that
$$\cases
2\partial_\tau w+(V-w)\p_q w=0,\qquad\text{$(q,\tau)\in (-\infty,M]\times[0,\tau_0)$},\\
w(q,0)=F_0'(q).\\
\endcases\tag 2.2
$$

The characteristics of (2.2) starting at the point $(0, s)$ is defined by
$$\cases
\ds\f {dq}{d\tau}(\tau,s)=\f{1}{2}(V-w)(q(\tau,s),\tau),\\
q(0,s)=s.\\
\endcases\tag2.3
$$

Along this characteristic curve, we have
$$\cases
\ds\frac{dw}{d\tau}(q(\tau,s),\tau)=0,\\
w(q(0,s), 0)=F_0'(s),\\
\endcases
$$
which yields for $\tau<\tau_0$
$$w(q(\tau,s),\tau)=F_0'(s)=\partial_q V(q(\tau,s),\tau).\tag2.4$$

>From (2.3)-(2.4), we obtain
$$\cases
\partial_{\tau s}q(\tau,s)=\ds\f{1}{2}F_0'(s)\p_s q(\tau,s)-\f{1}{2}F_0''(s),\\
\partial_s q(0,s)=1.\\
\endcases
$$
This yields $\p_sq(\tau,s)=\ds\exp\biggl(\ds\f{1}{2}F_0'(s)\tau\biggr)
\biggl(1-\f{F_0''(s)}{F_0'(s)}\biggr)
+\f{F_0''(s)}{F_0'(s)}>0$ if $F_0'(s)\not=0$ and
$\p_sq(\tau,s)=1-\ds\f12\tau F_0''(s)>0$ if $F_0'(s)=0$ when
$0\leq \tau<\tau_0$, and then
$$
q(\tau,s)=q(\tau,M)+\int_M^s\biggl(\ds \exp\biggl(\ds\f{1}{2}F_0'(\rho)\tau\biggr)
\biggl(1-\f{F_0''(\rho)}{F_0'(\rho)}
\biggr)
+\f{F_0''(\rho)}{F_0'(\rho)}\biggr)d\rho,\tag2.5
$$
and
$$
V(q(\tau,s),\tau)=2\p_\tau q(\tau,s)+w=2\p_\tau q(\tau,M)
+\int_M^s \exp\biggl(\ds\f{1}{2}F_0'(\rho)\tau\biggr)(F_0'(\rho)-F_0''(\rho))
\,d\rho+F_0'(s),\tag2.6
$$
here we have used that $\ds\lim_{z\to
0}\bigl(e^{z\tau} \bigl(1-\f{y}{z}\bigr)+\f{y}{z}\bigr)=1-\tau y$.

Note that $q(\tau, M)=M$ such that $V(q,\tau)$ satisfies the boundary
condition $V|_{q=M}=0$.  Together with (2.5)-(2.6), this yields
$V(q(\tau,s),\tau)=\ds\int_M^s \exp\biggl(\ds\f{1}{2}F_0'(\rho)\tau\biggr)(F_0'(\rho)-F_0''(\rho))\,d\rho+F_0'(s)$
and
$q(\tau,s)=M+\ds\int_M^s\biggl(\ds \exp\biggl(\ds\f{1}{2}F_0'(\rho)\tau\biggr)
\biggl(1-\f{F_0''(\rho)}{F_0'(\rho)}\biggr)
+\f{F_0''(\rho)}{F_0'(\rho)}\biggr)\,d\rho$. On the other hand, by the
implicit function theorem, we can obtain the smooth function
$s=s(q,\tau)$ for $\tau<\tau_0$. Therefore, $V(q,\tau)=F_0(s(q,\tau))$
is a smooth solution of (2.1) for $0\leq \tau<\tau_0$.
\qed \enddemo

We now start to construct an approximate solution of (1.3) for
$0\leq\tau=\ve\ln (1+t)<\tau_0$.

Let $w_0$ be the solution of the linear wave equation
$$\cases
\p_t^2 w-\triangle w=0,\\
w(0,x)=u_0(x),\\
\p_t w(0,x)=u_1(x).\\
\endcases
$$

Choose a $C^\infty$ function $\chi(s)$ such that $\chi(s)=1$ for
$s\leq 1$ and $\chi(s)=0$ for $s\geq 2$.  We set, for
$0\leq \tau=\ve\ln (1+t)<\tau_0$,
$$
u_a(t,x)=\ve\chi(\ve t)w_0(t,x)+\f{\ve}{r}\left(1-\chi(\ve t)\right)
\chi(-3\ve q)V(q,\tau).\tag 2.7
$$
By [9,Theorem 6.2.1] and Lemma 2.1, we have $|Z^\al u_a|\leq
C_{\al,b}\ve \left(1+t\right)^{-1}$ for $\tau\leq b<\tau_0$ and all
multi-indices $\alpha$. Set
$$J_a=\p_t^2 u_a-(1+u_a+\p_tu_a)\triangle u_a.$$

\proclaim{Lemma 2.2}
$$
\int_0^{e^{b/\ve}-1}\|Z^\alpha J_a\|_{L^2}\,dt\leq C_{\al, b}\ve^{3/2}|\ln\ve|.
$$
\endproclaim

\demo{Proof} We divide the proof procedure into the following three cases.

{\it Case 1.} $\ds\f{2}{\varepsilon}\leq t\leq e^{b/\varepsilon}-1$.
In this case, $\chi(\ve t)=0$ and $u_a(t,x)=\ds\f{\ve}{r}\,\chi(-3\ve
q)V(q,\tau)$. Then
$$
J_a=-\f{\ve^2}{r^2}\bigl(\bar{V}\p_q^2\bar{V}-\p_q\bar{V}\p_q^2\bar{V}+2\p_{q\tau}
\bar{V}\bigr)+O\biggl(\f{\ve^2}{(1+t)^3}\biggr),
$$
where $\bar{V}(q,\tau)=\chi(-3\ve q)V(q,\tau)$. Since $\p_qV$ has
compact support, we have
$$
\multline
\quad(\bar{V}\p_q^2\bar{V}-\p_q\bar{V}\p_q^2\bar{V}+2\p_{q\tau} \bar{V})(q,\tau)\\
=9\ve^2\chi(-3\ve q)\chi''(-3\ve q)V^2(q,\tau)+27\ve^3\chi'(-3\ve q)\chi''(-3\ve q)V^2(q,\tau)
-6\ve\chi'(-3\ve q)\p_\tau V(q,\tau).
\endmultline
$$
Hence $|Z^\al J_a|\leq C_{\al,b}\ve^2\left(1+t\right)^{-3}+C_{\al,b}\ve^3
\left(1+t\right)^{-2}\psi(-3\ve q)$, where
$\psi(s)$ is a cut-off function satisfying $\psi(s)=1$ for $s\in [1,
2]$, and $\psi(s)=0$ otherwise.

{\it Case 2.} $t\leq \ds\f{1}{\ve}$.  In this case, $\chi(\ve t)=1$
and $u_a=\ve w_0$. This gives $J_a=-\ve^2(w_0+\p_tw_0)\triangle
w_0$. It follows then from a direct computation that
$$
|Z^\alpha J_a|\leq C_\alpha\ve^2\left(1+t\right)^{-2}.
$$

{\it Case 3.} $\ds\f{1}{\varepsilon}\leq t\leq\ds\f{2}{\varepsilon}$.
A direct computation yields
$$
u_a=\ve w_0+\ve \left(1-\chi(\ve t)\right)\left(r^{-1}\chi(-3\ve q)V-w_0\right).
$$
Then
$$
J_a=J_1+J_2+J_3+J_4,
$$
where
$$\align
J_1&=-(u_a+\p_tu_a)\triangle u_a,\\
J_2&=\ve(\p_t^2-\triangle)\bigl[(1-\chi(\ve t))r^{-1}\chi(-3\ve q)(V-F_0)\bigr],\\
J_3&=\ve(\p_t^2-\triangle)\bigl[(1-\chi(\ve t))\chi(-3\ve q)(r^{-1}F_0-w_0)\bigr],\\
J_4&=\ve(\p_t^2-\triangle)\bigl[(1-\chi(\ve t))\left(\chi(-3\ve q)-1\right)w_0\bigr].
\endalign
$$

It is easy to see $|Z^\alpha J_1|\leq C_{\alpha,b}\ve^2\left(1+t\right)^{-2}$.

Due to
$(\p_t^2-\p_r^2)(V(q,\tau)-F_0(q))=(\p_t-\p_r)(\p_{\tau}V\ds\f{\ve}{1+t})
=-2\p_{\tau
q}^2V\ds\f{\ve}{1+t}+\p_{\tau}^2V\ds(\f{\ve}{1+t})^2-\p_{\tau}V\f{\ve}{(1+t)^2}$,
and $V(q,0)=F_0(q)$, we have $|Z^\alpha J_2|\leq
C_{\alpha,b}\ve^2\,|\ln\ve|\left(1+t\right)^{-2}$.

Moreover, by [9, Theorem 6.2.1], we have that, for any constant $l>0$,
if $r\geq lt$, then
$$
|Z^\alpha(w_0-r^{-1}F_0)|\leq C\left(1+t\right)^{-2}.\tag2.8
$$
On the other hand, from $\p_t=\ds\f{tS-rH}{t^2-r^2}$ and
$\p_r=\ds\f{tH-rS}{t^2-r^2}$ we obtain
$\triangle=\ds\f{1}{r+t}(S+H)(\p_t-\p_r)-\f{2}{r}\p_r$.  Therefore,
$|Z^\alpha J_3|\leq C_{\alpha,b}\ve\left(1+t\right)^{-3}\leq
C_{\alpha,b}\ve^2\left(1+t\right)^{-2}$.

Since the support of $J_4$ is in ${q\leq-
\ds\f{1}{3\ve}}$, combining the fact that, for any $\phi(t,r)\in C^1$,
$$|\p \phi|\le\ds\f{C}{1+|t-r|}\dsize\sum_{|\beta|=1}|Z^{\beta}\phi|,\tag2.9$$
we get the estimate $|Z^\alpha J_4|\leq C_{\alpha,b}\ve^3\left(1+t\right)^{-1}
\leq C_{\alpha,b}\ve^2\left(1+t\right)^{-2}$.

\smallskip

The above analysis yields
$$
|Z^\al J_a|\leq C_\al\ve^2\left(1+t\right)^{-2}|\ln\ve|.
$$
Collecting the estimates above, we arrive at
$$
\aligned
\|Z^\al J_a\|_{L^2}&\leq C_{\al, b}\ve^{5/2}\left(1+t\right)^{-1}, \qquad\enspace
 \f{2}{\ve}\leq t\leq e^{b/\ve}-1,\\
\|Z^\al J_a\|_{L^2}&\leq C_\al\ve^2\,|\ln\ve|\left(1+t\right)^{-1/2}, \quad
t\leq\f{2}{\ve}.
\endaligned
$$
Consequently,
$$
\int_0^{e^{b/\ve}-1}\|Z^\al u_a\|_{L^2}\,dt\leq C_{\al,b}\ve^{3/2}\,|\ln\ve|,
$$
and Lemma 2.2 is proved. \qed
\enddemo

For later reference, we cite a result that was shown in [17].

\proclaim{Lemma 2.3} Let $f(t,x)\in C^1(\Bbb R^4)$ only depend on $(t, r)$. Moreover,
assume that $\operatorname{supp} f\subseteq\{(t, x): r\leq M+t\}$. Then
$$
\|(1+|t-r|)^{-1}f\|_{L^2}\leq C\|\p_r f\|_{L^2}.
$$
\endproclaim

Based on these preparations, we next establish:

\proclaim{Proposition 2.4} For sufficiently small $\varepsilon>0$ and
$0\leq \tau=\ve\ln(1+t)\leq b<\tau_0$,
\rom{(1.3)} has a $C^\infty$ solution $u(t,x)$ which satisfies, for all $|\alpha|\le
2$,
$$
|Z^\alpha\p (u-u_a)|\leq
C_{b}\varepsilon^{3/2}\,|\ln\varepsilon|(1+t)^{-1}(1+|t-r|)^{-1/2}.\tag2.10
$$
\endproclaim

\demo{Proof} Set $v=u-u_a$. Then one has
$$\cases
\p_t^2 v-(1+u+\p_tu)\triangle v=-J_a+(v+\p_tv)\triangle u_a,\\
v(0,x)=\p_t v(0,x)=0.
\endcases\tag 2.11
$$

We make the induction hypothesis that, for some  $T\leq e^{b/\ve}-1$,
$$
|Z^\alpha\p v|\leq\varepsilon(1+t)^{-1}(1+|t-r|)^{-1/2},\quad |\alpha|\leq 2,
\quad t\leq T, \tag 2.12
$$
which then further implies that, for $|\alpha|\leq 2$ and $t\leq T$,
$$|Z^\alpha v|\leq C\varepsilon(1+t)^{-1}(1+|t-r|)^{1/2}.\tag2.13
$$

To verify the validity of (2.12), we will prove that, for sufficiently small $\varepsilon>0$,
$$
|Z^\alpha\p
v|\leq\f{\varepsilon}{2}\left(1+t\right)^{-1}(1+|t-r|)^{-1/2},\quad
|\alpha|\leq 2, \quad t\leq T,\tag2.14
$$
and apply the continuity method to obtain $\ve\ln(1+T)=b$.

Applying  $Z^\alpha$ to both sides of (2.11) and using $\left[Z^{\al},
\p_t^2-\triangle\right] =\dsize\sum_{|\beta|<|\al|}C_{\al\beta}
Z^{\beta}\left(\p_t^2-\triangle\right)$ yields, for any $|\alpha|\leq 4$,
$$
\left(\p_t^2-\triangle\right)Z^\alpha v
=Z^\alpha G-\sum_{|\beta|<|\alpha|}C_{\alpha\beta}Z^\beta G,\tag2.15
$$
where
$$
G=\left(u+\p_tu\right)\triangle v-J_a+\left(v+\p_tv\right)\triangle u_a.
$$
Since
$$
Z^\al G=\left(u+\p_tu\right)Z^\al\triangle v
+\sum_{\Sb\al_1+\al_2=\al\\|\al_1|\geq 1\endSb}Z^{\al_1}(u+\p_tu)Z^{\al_2}
\triangle v+Z^\al\left[-J_a+(v+\p_tv)\triangle u_a\right],
$$
we have from (2.15) that
$$
\left(\p_t^2-(1+u+\p_tu)\triangle\right)Z^\al v=F,\tag2.16
$$
where
$$
\align
F=&\sum_{\Sb\al_1+\al_2=\al\\|\al_1|\geq 1\endSb}Z^{\al_1}(u+\p_tu)Z^{\al_2}\triangle v
+Z^\al\left[-J_a+(v+\p_tv)\triangle u_a\right]
-\sum_{|\beta|<|\alpha|}C_{\alpha\beta}Z^\beta G+(u+\p_tu)\left[Z^\alpha,\triangle\right]v.
\endalign
$$

Next we derive an estimate of $\|\p Z^\alpha v\|_{L^2}$ from equation (2.16).
Define the energy
$$
E(t)=\f{1}{2}\sum_{|\al|\leq 4}\int_{\Bbb R^3}
\left(|\p_t Z^\alpha v|^2+(1+u+\p_tu)|\nabla Z^\al v|^2\right)dx.
$$
Multiplying both sides of (2.16) by $\p_t Z^\alpha v$ ($|\al|\le 4$)
and integrating by parts, and noting that $|\p^\beta u|=|\p^\beta
u_a+\p^\beta v|\leq C_b\ve\left(1+t\right)^{-1}$ ($|\beta|=1,2$), which
follows from the construction of $u_{a}$ and assumption (2.12), we
arrive at
$$
E'(t)\leq \f{C_b\ve}{1+t}\,E(t)
+\sum_{|\al|\leq 4}\int_{\Bbb R^3}\left|F\right|\left|\p_t Z^\al v\right|dx.\tag2.17
$$

We now treat each term arising in the integration of
$\dsize\sum_{|\al|\leq 4}\int_{\Bbb R^3}\left|F\right|\left|\p_t Z^\al
v\right|dx$ separately.

(A) {\it The term $\dsize\sum_{\Sb\al_1+\al_2=\al\\|\al_1|\geq
1\endSb} \int_{\Bbb R^3}\left|Z^{\al_1}(u+\p_tu)Z^{\al_2}\triangle v\right|
\left|\p_t Z^\alpha v\right|dx$.}
First, we note that:
\roster
\item"(i)" By assumption (2.13), for $|\beta|\le 2$, we have
$$\left|(1+|t-r|)^{-1}Z^{\beta}(u_a+v)\right|\le \ds\f{C_b\ve}{1+t}.\tag2.18$$

\item"(ii)" By (2.18) and (2.9), for $|\al_1|+|\al_2|=|\al|\leq 4$ with
$|\al_1|\ge 1$, we have
$$
\align
& \int_{\Bbb R^3}\left|Z^{\al_1}(u+\p_tu)(Z^{\al_2}\triangle v)\p_t
Z^\alpha v\right|dx\\
& \enspace \leq C_b\int_{\Bbb
R^3}\left|(Z^{\al_1}u_a+Z^{\al_1}\p_tu_a)(Z^{\al_2}\triangle v)\p_t Z^\al
v\right|dx+\int_{\Bbb R^3}\left|(Z^{\al_1}v+Z^{\al_1}\p_tv)(Z^{\al_2}\triangle
v)\p_t Z^\al v\right|dx\\
& \enspace \leq\f{C_{b}\ve}{1+t}\,E(t)+C\sum_{|\g|\leq
|\al_2|+1}\int_{\Bbb R^3}\left|(1+|t-r|)^{-1}Z^{\al_1} v\right|\left|Z^\g\p
v\right|\left|\p_t Z^\al v\right|dx
\\ &\qquad+\,C_b\sum_{|\g|\leq|\al_2|+1}\int_{\Bbb
R^3}\left|Z^{\al_1} \p_tv\right|\left|Z^\g\p v\right|\left|\p_t Z^\al v\right|dx.\tag2.20
\endalign
$$
\endroster
Note further that there is at most one number larger than 2 between
$|\al_1|$ and $|\g|$.  If $|\al_1|>2$, then $|\g|\le 2$ and, by
Lemma 2.3 applied to $(1+|t-r|)^{-1}Z^{\al_1} v$ and by assumption
(2.12), we arrive at
$$
\int_{\Bbb R^3}
\left|(1+|t-r|)^{-1}Z^{\al_1} v\right|\left|Z^\g\p v\right|\left|\p_t Z^\al v\right|dx
+\int_{\Bbb R^3}\left|Z^{\al_1} \p_tv\right|\left|Z^\g\p v\right|\left|\p_t Z^\al
v\right|dx \leq \f{C_b\ve}{1+t}\,E(t).\tag2.21
$$

If $|\g|>2$, then $|\al_1|\le 2$ and it follows from (2.13) that
$\left|(1+|t-r|)^{-1}Z^{\al_1} v\right|\leq
C\ve(1+t)^{-1}(1+|t-r|)^{-1/2}$ which leads to
$$
\int_{\Bbb R^3}\left|(1+|t-r|)^{-1}Z^{\al_1} v\right|\left|Z^\g\p v\right|
\left|\p_t Z^\al v\right|dx
+\int_{\Bbb R^3}\left|Z^{\al_1} \p_tv\right|\left|Z^\g\p v\right|\left|\p_t Z^\al
v\right|dx \leq \f{C_b\ve}{1+t}\,E(t).\tag2.22
$$

Inserting (2.21)--(2.22) into (2.20) yields
$$\dsize\sum_{\Sb\al_1+\al_2=\al\\|\al_1|\geq 1\endSb}\int_{\Bbb R^3}
\left|Z^{\al_1}(u+\p_tu)Z^{\al_2}\triangle v\right|\left|\p_t Z^\alpha v\right|dx
\le\ds\f{C_b\ve}{1+t}\,E(t).\tag2.23$$

(B) {\it The terms $\int_{\Bbb R^3}|Z^\beta\bigl((u+\p_tu)\triangle
v\bigr)\cdot\p_tZ^\al v|dx$ with $|\beta|<|\alpha|$.} We only need to
treat the term $\int_{\Bbb R^3}\left|(u+\p_tu)Z^\beta\triangle
v\cdot\p_tZ^\al v\right|dx$, since the other terms have been estimated in
(A).

By (2.12), we have
$$
\align
&\int_{\Bbb R^3}\left|(u+\p_tu)Z^\beta\triangle v\cdot\p_tZ^\al v\right|dx\\
&\qquad
\leq C\sum_{|\g|\leq|\beta|+1\leq|\al|}\int_{\Bbb R^3}\left|(1+|t-r|)^{-1}u\right|
\left|Z^\g\p v\right|\left|\p_tZ^\al v\right|dx
+\int_{\Bbb R^3}\left|\p_tu\right|\left|Z^\beta\triangle v\right|
\left|\p_tZ^\al v\right|dx\\
&\qquad \leq\f{C_{b}\ve}{1+t}\,E(t).\tag2.24
\endalign
$$

(C) {\it The terms
$\int_{\Bbb R^3}\left|\p_t Z^\al v\right|\left|Z^\beta J_a\right|dx$
with $|\beta|\leq|\al|\leq 4$.}
In this case, we have
$$
\int_{\Bbb R^3}\left|\p_t Z^\al v\right|\left|Z^\beta J_a\right|dx
\leq\|Z^\beta J_a\|_{L^2}\sqrt{E(t)}.\tag2.25
$$

(D) {\it The terms
$\int_{\Bbb R^3}\left|Z^\beta((v+\p_tv)\triangle u_a)\right|\left|\p_t Z^\al v\right|dx$
with $|\beta|\leq|\al|\leq 4$.} A direct computation yields
$$
\align
&\int_{\Bbb R^3}\left|Z^\beta((v+\p_tv)\triangle u_a)\right|\left|\p_t Z^\al v\right|dx\\
&\qquad \leq C\sum_{\Sb|\beta_1|+|\beta_2|\leq|\beta|+1\\|\beta_1|\leq|\beta
|\endSb}\int_{\Bbb R^3}\big(\left|(1+|t-r|)^{-1}(Z^{\beta_1}v)(Z^{\beta_2}\p u_a)\right|
+\left|(Z^{\beta_1}\p_tv)(Z^{\beta_2}\p u_a)\right|\big)\left|\p_t Z^\al v\right|dx\\
&\qquad \leq\f{C_{b}\ve}{1+t}\,E(t).\tag2.26
\endalign
$$

(E) {\it The term $\int_{\Bbb
R^3}\left|(u+\p_tu)\left[Z^\al,\triangle\right]v\right|\left|\p_t
Z^\al v\right|dx$.}  Since
$\left[Z^\al,\triangle\right]=\ds\sum_{|\beta|\leq|\al|-1}C_{\al\beta}\p^2Z^\beta$,
we have
$$
\align
\int_{\Bbb R^3}\left|(u+\p_tu)\left[Z^\al,\triangle\right]v\right|
\left|\p_t Z^\al v\right|dx
& \leq C\sum_{|\beta|\leq|\al|-1}\int_{\Bbb R^3}\left|u+\p_tu\right|
\left|\p^2Z^\beta v\right|\left|\p_t Z^\al v\right|dx\\
& \leq C\sum_{\Sb|\beta|\leq|\al|-1\\|\g|\leq|\beta|+1\endSb}\int_{\Bbb R^3}
\left|(1+|t-r|)^{-1}
(u+\p_tu)\right|\left|\p Z^\g v\right|\left|\p_tZ^\al v\right|dx\\
& \leq\f{C_{b}\ve}{1+t}\,E(t).\tag2.27
\endalign
$$

\smallskip

Substituting (2.23)--(2.27) into (2.17) yields
$$
E'(t)\leq\f{C_{b}\ve}{1+t}\,E(t)
+\ds\sum_{|\beta|\leq 4}\|Z^\beta J_a\|_{L^2}\sqrt{E(t)}.\tag2.28
$$
Thus, by Lemma 2.2 and Gronwall's inequality we obtain
$$
\|\p Z^\al v\|_{L^2}\leq C_{b}\ve^{3/2}\,|\ln\ve|,\quad|\al|\leq 4,
$$
and then
$$
\|Z^\al\p v\|_{L^2}\leq C_{b}\ve^{3/2}\,|\ln\ve|,\quad|\al|\leq 4.\tag2.29
$$

By (2.29) and the Klainerman-Sobolev inequality (see [9] or [14]), we have
$$
|Z^\al\p v|\leq C_{b}\ve^{3/2}\,|\ln\ve|(1+t)^{-1}(1+|t-r|)^{-1/2},\quad |\alpha|\leq 2,
\enspace t\leq T.\tag2.30
$$
This implies that, for $\ve>0$ small enough,
$$
|Z^\alpha\p v|\leq\f{\varepsilon}{2}\,(1+t)^{-1}(1+|t-r|)^{-1/2},\quad |\alpha|\leq 2,
\enspace t\leq T.
$$

Therefore, we have completed the proof of (2.12) and, together with
(2.30), the proof of (2.10).
\qed
\enddemo

Proposition 2.4 immediately gives $\ds\lim_{\overline{\ve\rightarrow
0}}\ve\ln(1+ T_\ve)\geq\tau _0$ and, therefore,
$$\ds\lim_{\overline{\ve\rightarrow 0}}\ve\ln T_\ve\geq\tau
_0.\tag2.31$$

\vskip 0.3 true cm
\remark{Remark 2.1}
The analysis of this section can be directly applied to problem
(1.6) with general initial data so that a lower bound on the
lifespan $T_{\ve}$ is obtained as in [5].
\endremark


\head \S 3. The upper bound on the
lifespan $T_\ve$ and proof of Theorem 1.1. \endhead

In this section, we focus on the upper bound on $T_{\ve}$. Here some
ideas are inspired by [10] and [6]. Since equation (1.3) contains the
solution $u$ and its derivatives simultaneously rather than the
derivatives of $u$ only, as in [10], and the function $u$ only, as in
[6], respectively, we have to be more careful in computations and also
need to estimate more quantities. Thanks to the estimate of
$Z^{\al}(u-u_a)$ with $|\al|\le 2$ in (2.13), we observe that
$|Z^{\al}(u-u_a)|\le C\varepsilon\left(1+t\right)^{-1}$ holds near the
light cone which plays an important role in the analysis below.

Let $U=ru$ and $c^2(u)=1+u+\p_tu$. Due to radial symmetry of $u$,
(1.3) can be rewritten as
$$\cases
\enspace \partial_t^2 U-c^2\p_r^2U=0, \\
\enspace U(0, r)=\ve ru_0,\quad
\p_t U(0, r)=\ve ru_1.\\
\endcases\tag3.1
$$
Define two operators
$$
L_1=\p_t+c\p_r,\q L_2=\p_t-c\p_r.
$$
We also set
$$
w_1=L_2\p_rU=\left(\p_t-c\p_r\right)\p_rU,\q
w_2=L_1\p_rU=\left(\p_t+c\p_r\right)\p_rU,
$$
which yields $\p_{tr}U=\ds\f{w_1+w_2}{2}$ and $\p^2_{r}U=\ds\f{w_2-w_1}{2c}$.

Note that
$$
L_1L_2=\p_t^2-c^2\p_r^2-(L_1c)\p_r,\q L_2L_1=\p_t^2-c^2\p_r^2+(L_2c)\p_r.
$$
Then
$$
\align
L_1 w_1&=L_1L_2\p_rU
=-\f{1}{4rc}w_1^2+\f{w_1}{4rc}(w_2+\p_tu+\f{r}{c}L_2u)-\f{w_2}{4rc}
(\p_tu+\f{r}{c}L_2u),\tag 3.2\\
L_2 w_2&=L_2L_1\p_rU
=\f{1}{4rc}w_2^2+\f{w_2}{4rc}(-w_1-\p_tu+\f{r}{c}L_1u)+\f{w_1}{4rc}
(\p_tu-\f{r}{c}L_1u).\tag 3.3
\endalign
$$
Due to $\p_r
c=c'\p_ru=\ds\f{1}{2c}\p_ru-\f{1}{2cr}\p_tu+\f{1}{4rc}(w_1+w_2)$, we
have
$$
\align
L_1w_1+w_1\p_r c&=\f{w_1}{4rc}(2w_2-\p_tu+\f{r}{c}L_1u)-\f{w_2}{4rc}(\p_tu+\f{r}{c}L_2u),\\
L_2w_2-w_2\p_r c&=\f{w_2}{4rc}(-2w_1+\p_tu+\f{r}{c}L_2u)+\f{w_1}{4rc}(\p_tu-\f{r}{c}L_1u)
\endalign
$$
and
$$
\align
d(|w_1|(dr-cdt))&=\operatorname{sgn}w_1\left(L_1w_1+w_1\p_r c\right)dt\wedge dr\\
&=\operatorname{sgn}w_1\left[\f{w_1}{4rc}(2w_2-\p_tu+\f{r}{c}L_1u)
-\f{w_2}{4rc}(\p_tu+\f{r}{c}L_2u)\right]dt\wedge dr,\tag3.4 \\
d(|w_2|(dr+cdt))&=\operatorname{sgn}w_2\left(L_2w_2-w_2\p_r c\right)dt\wedge dr\\
&=\operatorname{sgn}w_2\left[\f{w_2}{4rc}(-2w_1+\p_tu+\f{r}{c}L_2u)
+\f{w_1}{4rc}(\p_tu-\f{r}{c}L_1u)\right]dt\wedge dr.\tag 3.5
\endalign
$$

$$\epsfysize=99mm \epsfbox{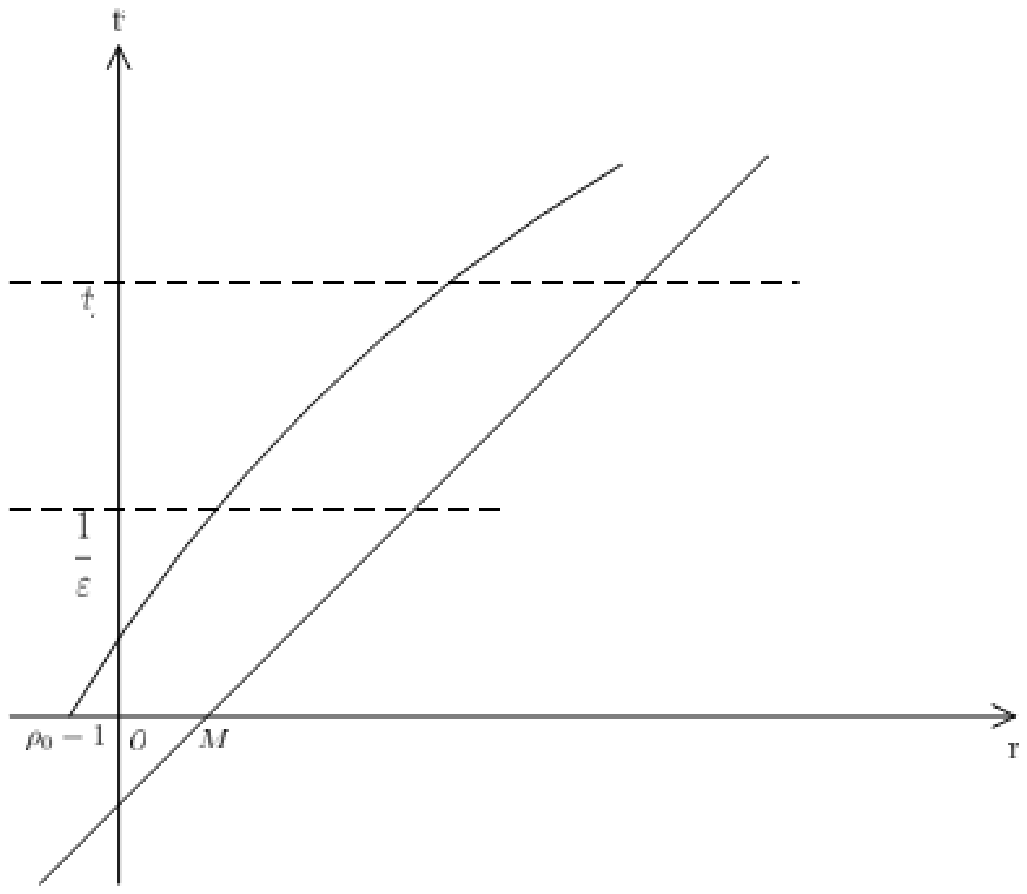}$$
\centerline{\bf Figure 1.}

\vskip 0.6 true cm

>From $\S 2$ we have that when $\ve\ln (1+T_b)=b<\tau_0$ and $b>0$ is a
fixed constant, (1.3) has a $C^\infty$ solution for $t\leq
T_b$. Choosing $\ve$ sufficiently small such that $1/\ve<e^{b/\ve}-1$,
we define the characteristics $\G_\la^\pm$ by $dr/dt=\pm c$,
respectively, passing through $(\la,0)$ in the $(r,t)$ plane. Let
$D$ be the domain which is bounded by $\G_M^+$ and $\G_{\rho_0-1}^+$
(see Figure 1), where $\rho_0$ is chosen such that
$\tau(\rho_0)=\ds\f{2}{F_0'(\rho_0)}
\ln\f{F_0''(\rho_0)}{F_0''(\rho_0)-F_0'(\rho_0)}=\tau_0$
or $\ds\f{2}{F_0''(\rho_0)}=\tau_0$. (Another possibility is that
there exists no point $\rho_0\in (-M, M)$ such that
$\tau(\rho_0)=\tau_0$. In this case, by a minor modification which
uses the fact that, for any fixed $\dl_0>0$ sufficiently small, there
exists $\rho_0^{\dl}\in (-M, M)$ such that $\tau
(\rho_0^{\dl})\le\tau_0+\dl_0$, an analogous proof still works).
Obviously, $\G_M^+$ is the straight line $r=t+M$.

Fix a positive constant $R$ satisfying $R>\tau_0$, and assume that $u$
is a solution of (1.3) (at least) for $t\leq T\leq \exp(R/\ve)$. We
define
$$
\align
 A(t)&=\sup_{1/\ve\leq s\leq t}\int_{(s,r)\in D}|w_1(s,r)|\,dr,\\
 B(t)&=\sup_{\Sb 1/\ve\leq s\leq t\\(s,r)\in D\endSb}s\left|\p_su(s,r)\right|,\\
 C(t)&=\sup_{\Sb 1/\ve\leq s\leq t\\(s,r)\in D\endSb}s^{3/2}\left|w_2(s,r)\right|,\\
 D(t)&=\sup_{\Sb 1/\ve\leq s\leq t\\(s,r)\in D\endSb}s^{3/2}\left|L_1u(s,r)\right|.
\endalign
$$

\proclaim{Lemma 3.1}
There exists a constant $E>0$ such that, for $\ve>0$ sufficiently small,
$$A(1/\ve)\leq \ds\f{E\ve}{2},\quad B(1/\ve)\leq E\ve,\quad C(1/\ve)\leq E^2\ve^2,
\quad D(1/\ve)\leq E\ve.\tag3.6$$
\endproclaim

\demo{Proof} The equation $r=r(t)$ of $\G^+_{\la}$ for $\la\in
[\rho_0-1,M]$ is
$$
\cases
&\ds\f{dr(t)}{dt}=c(u(t,r(t)))\equiv c(t),\\
&r(0)=\la.
\endcases
$$
This gives
$$r(t)-\la=\int_0^t(c(s)-1)\,ds+t.$$ Because of $|c(t)-1|\le
C_b\left(|u(t,r(t))|+|\p_tu(t,r(t))|\right)\le C_b\ve
(1+t)^{-1}(1+|t-r(t)|)^{1/2}$ for $0\leq \tau=\ve\ln(1+t)\leq
b<\tau_0$,
$$
\align
|r(t)-t|&\le |\la|+C_b\int_0^t\ve \left(1+s\right)^{-1}(1+|s-r(s)|)^{1/2}\,ds\\
&\le m_0+b+C_b\int_0^t\ve \left(1+s\right)^{-1}|s-r(s)|\,ds,
\tag3.7
\endalign
$$
where $m_0=\max\left\{|\rho_0-1|,M\right\}$. It follows from
Gronwall's inequality that
$$|r(t)-t|\le C_b.$$
Therefore,
$$\left|t+M-r(t)\right|\le C_b,\tag3.8$$
which means the distance between $\G_{\rho_0-1}^+$ and $\G_M^+$ is bounded.

The equation $\widetilde{r}=\widetilde{r}(t)$ of $\G^-_{\mu}$ is
$$
\cases
&\ds\f{d\widetilde{r}(t)}{dt}=-c(u(t,\widetilde{r}(t)))\equiv -c(t),\\
&\widetilde{r}(0)=\mu,
\endcases
$$
and then $\left|\widetilde{r}(t)+t-\mu\right|\le C_b$ can analogously be shown.

For $(t,r),(t',r')\in\G^-_\mu\cap D$ for $\mu\in\Bbb R$ and
$(t,r)\in\G^+_{\la}$, $(t',r')\in\G^+_{\la'}$, we obtain
$$
|t-t'|\leq\f{1}{2}\left[|t+r-\mu|+|t'+r'-\mu|+|t-r-\la|
+|t'-r'-\la'|+|\la-\la'|\right]
\leq C_b.\tag3.9
$$

We now start to prove (3.6).
Since $w_1=\p_tu+r\p_{tr}u-c\p_ru-cr\p_r^2u$, one has, for $t\le e^{b/\ve}-1$,
$$|w_1(t,r)|\le C_b\ve\left(1+|t-r|\right)^{-3/2}.\tag3.10$$
By (3.8)--(3.10), we thus arrive at
$$\int_{(1/\ve, r)\in D}|w_1(s,r)|\,dr\le C_b\ve.\tag3.11$$
In addition, because of $|\p_tu(t,r)|\le C_b\ve\left(1+t\right)^{-1}(1+|t-r|)^{-1/2}$
for $t\le e^{b/\ve-1}$, we have, for $(1/\ve, r)\in D$,
$$\f{1}{\ve}\left|\p_tu\biggl(\f{1}{\ve}, r\biggr)\right|\le C_b\ve.\tag3.12$$

Note that
$$w_2(t,r)=(\p_t+c\p_r)u+\f{r}{r+t}(S+H)\p_ru+(c-1)r\p_r^2u$$
which implies $|w_2(t,r)|\le C_b\ve(1+t)^{-1}$ for $(t,r)\in D$ and $\ve\ln (1+t)\le b$
in view of (2.13). Together with (3.3), this yields
$$|L_2w_2|\le\ds\f{C_b\ve^2}{(1+t)^2}.\tag3.13$$

Because of $w_2(t, t+M)=0$ and (3.9), we obtain from (3.13) that
$$C(1/\ve)\le C_b\ve^{5/2}\le C_b\ve^2.\tag3.14$$

Finally, from $L_1u=\ds\f{1}{r+t}(S+H)u+(c-1)\p_ru$ we have that
$|L_1u(t,r)|\leq C_b\ve(1+t)^{-2}$ for $t\le e^{b/\ve}-1$ and, therefore,
$$
D(1/\ve)\leq C_b\ve^{3/2}\leq C_b\ve.
$$

Collecting (3.11)--(3.12) and (3.14) completes the proof of
(3.6). \qed
\enddemo

Based on Lemma 3.1, we will use the continuity method to establish the
upper bound on $T_{\ve}$. To this end, we assume that, for $0\leq t\leq
T'\leq T$,
$$
A(t)\leq E\ve,\q B(t)\leq 2E\ve,\q C(t)\leq 3E^2\ve^2,\q D(t)\leq 2E\ve.\tag 3.15
$$

\pagebreak
We now start to verify:

\proclaim{Lemma 3.2} Under the assumption \rom{(3.15)} and if $\ve>0$ is
sufficiently small, then we have, for $1/\ve\leq t\leq T'$,
$$
A(t)\leq \f{2}{3}E\ve,\q B(t)\leq E\ve,\q C(t)\leq
\f{5}{2}E^2\ve^2,\q D(t)\leq E\ve.\tag3.16
$$
\endproclaim

\demo{Proof}
We first estimate $u(t,r)$ for $1/\ve\leq t\leq T'$ and $(t,r)\in
D$. Let $(t,r)\in\G^+_{\la}$. We integrate $L_1u$ along $\G^+_{\la}$
from time $1/\ve$ to $t$. From assumption (3.15), it is readily seen
that $|u(t,r)|\le C\ve^{3/2}$. Hence the definition
$c(t,r)=\sqrt{1+u+\p_tu}$ is possible and assumes a value close to
$1$. Therefore, $t\left|\p_ru(t,r)\right|\leq C\ve$ which yields
$|u(t,r)|\leq C\ve\left(1+t\right)^{-1}(1+|t-r|)$ for $1/\ve\leq t\leq
T'$ and $(t,r)\in D$.

Similar to the proof (3.8) and (3.9) we also have $|r-t|\le C$ and
$|t-t'|\le C$, where $t$ and $t'$ are same as in (3.9).

Now we estimate $A(t)$.  For $1/\ve\leq t\leq T'$, by equation (3.4)
and Green's formula, we have
$$
\align
&\int_{(t,r)\in D}|w_1(t,r)|\,dr\\
&\qquad\leq\int_{(1/\ve,r)\in D}|w_1(1/\ve,r)|dr+
\iint\limits_{\Sb 1/\ve\leq s\leq t\\(s,r)\in D\endSb}\left|\f{w_1}{4rc}
\left(2w_2-\p_tu+\f{r}{c}L_1u\right)-\f{w_2}{4rc}
\left(\p_tu+\f{r}{c}L_2u\right)\right|(s,r)\,dsdr\\
&\qquad\leq \f{E \ve}{2}+\iint\limits_{\Sb 1/\ve\leq s\leq t\\(s,r)\in
D\endSb}\left|\f{w_1}{4rc}\left(2w_2-\p_tu+\f{r}{c}L_1u\right)
-\f{w_2}{4rc}\left(\p_tu+\f{r}{c}L_2u\right)\right|(s,r)\,dsdr.\tag3.17
\endalign
$$
By the induction hypothesis (3.15), we have
$|\p_su(s,r)|\ds\leq\f{2E\ve}{s}$, $|w_2(s,r)|\ds\leq\f{3E^2\ve^2}{s^{3/2}}$,
and $|L_1u(s,r)|
\leq\ds\f{2E\ve}{s^{3/2}}$ for $1/\ve\le s\le T'$ and $(s,r)\in
D$. This immediately gives $|L_2u(s,r)|\leq\ds\f{5E\ve}{t}$. Note also
that $|r-s|\leq C$ holds for $s\geq 1/\ve$. We then have $r\geq s/2$
and
$$
\align
&\iint\limits_{\Sb 1/\ve\leq s\leq t\\(s,r)\in D\endSb}
\left|\f{w_2}{4rc}\left(\p_tu+\f{r}{c}L_2u\right)\right|(s,r)\,dsdr
\leq\iint\limits_{\Sb 1/\ve\leq s\leq t\\(s,r)\in D\endSb}\f{E^2\ve^2}{s^{3/2}r}
\left(\f{2E\ve}{s}+\f{6E\ve r}{s}\right)dsdr \\
& \qquad \leq25E^3\ve^3\int_{1/\ve}^t \f{ds}{s^{5/2}}\int_{(s,r)\in D}dr
\leq CE^3\ve^{9/2}\le \f{E\ve}{12},\tag3.18
\endalign
$$
where the generic constant $C>0$ is independent of $\ve$.

Similarly,
$$
\align
&\iint\limits_{\Sb 1/\ve\leq s\leq t\\(s,r)\in D\endSb}|\f{w_1}{4rc}(2w_2
-\p_tu+\f{r}{c}L_1u)|(s,r)dsdr
 \leq\iint\limits_{\Sb 1/\ve\leq s\leq t\\(s,r)\in D\endSb}\f{|w_1(s,r)|}{3r}
\left(\f{6E^2\ve^{2}}{s^{3/2}}+\f{2E\ve}{s}+\f{2E\ve r}{s^{3/2}c}\right)dsdr \\
&\qquad \leq4E\ve\int_{1/\ve}^t\f{ds}{s^{3/2}}
\int_{(s,r)\in D}|w_1(s,r)|{r}\,dr
\leq CE^2\ve^{5/2}.\tag3.19
\endalign
$$

Substituting (3.18)--(3.19) into (3.17) yields
$$
A(t)\leq\f{E\ve}{2}+\f{E\ve}{12}+CE^2\ve^{5/2},
$$
which implies that $A(t)\leq\f{2}{3}\,E\ve$ for $\ve>0$ sufficiently small.

Next we estimate $B(t)$.
Note that $u$ satisfies the equation
$$
L_2\p_tu=-\f{c}{r}\,w_1+\f{c}{r}\,\p_tu.\tag3.21
$$
$$\epsfysize=99mm \epsfbox{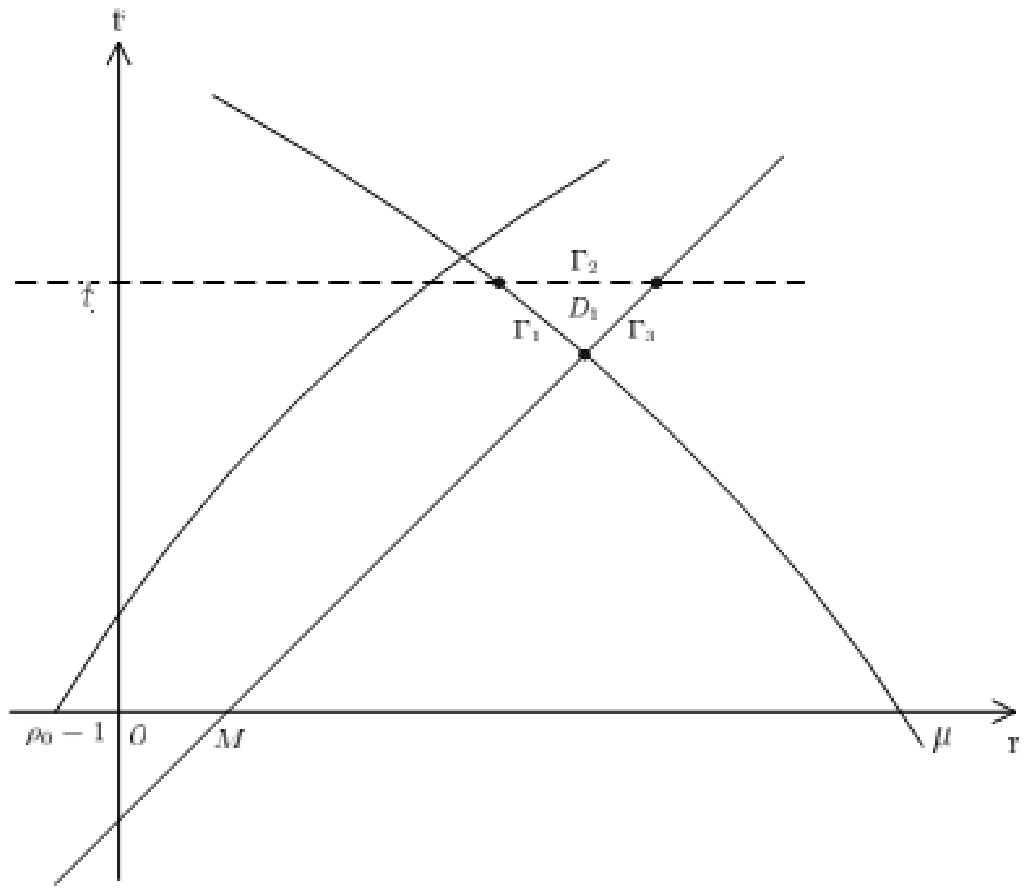}$$
\centerline{\bf Figure 2}

\vskip 0.6 true cm

We integrate equation (3.21) along the characteristics $\G_\mu^-$
which insects $\G_M^+$ at the point $(t',r')$.  For $t'\geq1/\ve$, we
denote by $D_1$ the domain bounded by $\G_\mu^-$, the line $\{t=t\}$,
and $\G_M^+$ (see the Figure 2) which are denoted by $\G_1$, $\G_2$
and $\G_3$, respectively. We have
$$
\iint\limits_{(s,r)\in D_1}\operatorname{sgn} w_1
\left[\f{w_1}{4rc}\left(2w_2-\p_tu+\f{r}{c}L_1u\right)
-\f{w_2}{4rc}\left(\p_tu+\f{r}{c}L_2u\right)\right](s,r)\,dsdr
=\left(\int_{\G_1}+
\int_{\G_2}\right)|w_1|\left(dr-cdt\right)
$$
which implies
$$
\align
& \int_{\G_1}|w_1|\left(dr-cds\right) \\
& \qquad \leq\iint\limits_{(s,r)\in D_1}\left|
\f{w_1}{4rc}\left(2w_2-\p_tu+\f{r}{c}L_1u\right)-\f{w_2}{4rc}\left(\p_tu
+\f{r}{c}L_2u\right)\right|(s,r)\,dsdr
+\int_{\G_2}|w_1|\left(dr-cdt\right)\\
& \qquad \leq \iint\limits_{\Sb 1/\ve\leq s\leq t\\(s,r)\in
D\endSb}\left|\f{w_1}{4rc}\left(2w_2-\p_tu+\f{r}{c}L_1u\right)
-\f{w_2}{4rc}\left(\p_tu+\f{r}{c}L_2u\right)\right|(s,r)\,dsdr
+\int_{(t,r)\in D_1}|w_1(t,r)|\,dr\\
& \qquad \leq\f{1}{6}\,E\ve+\f{2}{3}\,E\ve=\f{5}{6}\,E\ve.
\endalign
$$
This yields
$$
\align
\int_{t'}^t \f{|w_1(s,r(s))|}{r(s)}\,ds
&\leq\int_{t'}^t \f{|w_1(s,r(s))|}{r(t)}\,ds
\leq\f{2}{t}\int_{t'}^t|w_1(s,r(s))|\,ds\\
&=\f{2}{t}\int_{\G_1}\f{|w_1(s,r)|}{\sqrt{1+c^2}}\,ds
=\f{2}{t}\int_{\G_1}\f{|w_1|}{2c\sqrt{1+c^2}}\left(dr-cds\right)\\
&\leq\f{5}{6t}\,E\ve.
\endalign
$$
Note that $r(s)\geq\ds\f{s}{2}\geq\ds\f{1}{2\ve}$ holds.
>From (3.21) we then have
$$
\align
|\p_tu(t,r)|
&\leq\int_{t'}^t \f{|cw_1|(s,r(s))}{r(s)}ds
+\int_{t'}^t\f{\left|c\p_tu\right|(s,r(s))}{r(s)}\,ds\\
&\leq\f{11}{12t}\,E\ve+3\ve\int_{t'}^t |\p_tu(s,r(s)|\,ds.\tag3.22
\endalign
$$
By Gronwall's inequality, we obtain
$$
|u(t,r(t))|\leq\f{11}{12t}\,E\ve e^{C\ve}\leq \ds\f{E\ve}{t}.
$$
For $t'\leq1/\ve$, we have $t\leq t'+|t-t'|\leq 2/\ve$, and Section~2
tells us $|t\p_tu(t,r)|\leq E\ve$.  Therefore, $B(t)\le E\ve$ holds.

Finally, we estimate $C(t)$.
We write (3.3) as
$$L_2w_2=aw_2+b, \tag3.23$$
where
$$
a=\f{1}{4rc}\left(w_2-w_1-\p_tu+\f{r}{c}L_1u\right),\quad
b=\f{1}{4rc}w_1\left(\p_tu-\f{r}{c}L_1u\right).
$$
Integrating (3.3) along $\G_\mu^-$ as above, we have
$$
|w_2(t,r)|\leq\int_{t'}^t\left|aw_2+b\right|(s,r(s))\,ds.\tag3.24
$$
Noting $t'\geq t-|t-t'|\geq\ds\f{2t}{3}$, we obtain
$$
\align
\int_{t'}^t|b(s,r(s))|\,ds
&\leq\int_{t'}^t\left|\f{1}{4rc}w_1\p_tu\right|(s,r(s))\,ds
+\int_{t'}^t\left|\f{1}{4c^2}w_1L_1u\right|(s,r(s))\,ds\\
&\leq\f{2E\ve}{t^{1/2}}
\int_{t'}^t \f{|w_1(s,r(s))|}{r(s)}\,ds\leq\f{5E^2\ve^2}{3t^{3/2}}\tag3.25
\endalign
$$
and
$$
\int_{t'}^t\left|aw_2\right|(s,r(s))|\,ds\leq \f{E^2\ve^2}{t^{3/2}}\int_{t'}^t
\f{\left|w_2-w_1-\p_tu+\f{r}{c}L_1u\right|(s,r(s))}{r(s)}\,ds \leq
\f{CE^3\ve^{4}}{t^{3/2}}.\tag3.26
$$
Substituting (3.25)--(3.26) into (3.24) yields
$$
|w_2(t,r)|\leq\int_{t'}^t\left|aw_2+b\right|(s,r(s))\,ds\leq\f{5E^2\ve^2}{2t^{3/2}},
$$
and $C(t)\leq\ds\f{5E^2\ve^2}{2}$ is shown.

Eventually, we have that
$$
L_2L_1u=\f{2c^2}{r}\p_ru+\f{1}{2c}(L_2u)(\p_ru)
+\f{1}{2c}\left(-\f{c}{r}w_1+\f{c}{r}\p_tu\right)\p_ru.
$$
We integrate this equality along $\G_\mu^-$ from $t'$ to $t$ and obtain
$$\align
|L_1u(t,r)| &\leq\int_{t'}^t\f{6E\ve}{sr(s)}\,ds+\int_{t'}^t\f{15E^2\ve^2}{2s^2}\,ds
+\int_{t'}^t\left(\f{|w_1(s,r(s))|}{2r(s)}+\f{E\ve}{sr(s)}\right)\f{3E\ve}{s}\,ds\\
& \leq \f{CE\ve}{t^2}\leq\f{E\ve}{t^{3/2}}
\endalign
$$
and, therefore, $D(t)\leq E\ve$. \qed
\enddemo

Based on Lemmas 3.1--3.2, we will use the following blowup result to
establish the upper bound on the lifespan $T_{\ve}$.

\proclaim{Lemma 3.3}
Let w be a solution in $(0,T)$ of the ordinary differential equation
$$
\f{dw}{dt}=a_0(t)w^2+a_1(t)w+a_2(t)
$$
with $a_j$ continuous on $[0,T]$ and $a_0\geq 0$.
Let
$$
K=\left(\int_0^T|a_2(t)|\,dt\right)\exp\left(\int_0^T|a_1(t)\,dt\right).
$$
If $w(0)>K$, then
$$
\left(\int_0^Ta_0(t)\,dt\right)\exp\left(-\int_0^T|a_1(t)|\,dt\right)
<\left(w(0)-K\right)^{-1}.
$$
\endproclaim
\demo{Proof} The proof can be found in [8, Lemma~1.4.1] or [9, Lemma~1.3.2].
\qed \enddemo

Next we show that
$$ \overline{\lim_{\ve\rightarrow0}}\ve\ln
T_{\ve}\leq\tau_0.\tag3.27
$$

\demo{Proof of \rom{(3.27)}}
Along the characteristics $\G_{\rho_0}^+$, $w_1(t,r(t))$ satisfies
$$
\f{dw_1}{dt}(t,r(t))=L_1 w_1
=-\f{1}{4rc}w_1^2+\f{w_1}{4rc}\left(w_2+\p_tu+\f{r}{c}L_2u\right)-\f{w_2}{4rc}
\left(\p_tu+\f{r}{c}L_2u\right).
$$
We  fix $0<\mu<1$ and set $t_\ve=e^{\mu\tau_0/\ve}-1$ and $\widehat{w}(t,r(t))=-w_1(t,r(t))
\exp\biggl(-\ds\int_{t_\ve}^t\f{L_2u}{4c^2}(s,r(s))\,ds\biggr)$. Then
$$
\f{d\widehat{w}}{dt}(t,r(t))=a_0(t)\widehat{w}^2+a_1(t)\widehat{w}+a_2(t),
$$
where
$$\align
a_0(t)&=\f{1}{4r(t)c(t)}\exp\biggl(\ds\int_{t_\ve}^t\f{L_2u}{4c^2}(s,r(s))\,ds\biggr),\\
a_1(t)&=\f{w_2+\p_tu}{4rc}(t,r(t)),\\
a_2(t)&=-\biggl(\f{w_2}{4rc}
\left(\p_tu+\f{r}{c}L_2u\right)\exp\biggl(-\ds\int_{t_\ve}^t\f{L_2u}{4c^2}(s,r(s))
\,ds\biggr)\biggr)(t,r(t)).
\endalign
$$
By (3.15), we have that, for $0\leq t\leq T$,
$$
|a_1|\leq\f{E\ve}{t^2},\q |a_2|\leq\f{CE^3\ve^3}{t^{5/2}}
$$
which yields
$$
\int_{t_\ve}^T|a_1|ds\leq E\ve^2, \q
\int_{t_\ve}^T|a_2|ds\leq CE^2\ve^{9/2}.
$$
This further yields
$$
K=\left(\int_{t_\ve}^T|a_2(t)|\,dt\right)
\exp\left(\int_{t_\ve}^T|a_1(t)\,dt\right)=O(\ve^{9/2}).\tag3.28
$$

Due to $t_\ve\geq\ds 2/\ve$ and $ |3\ve(t_\ve-r(t_\ve)|\leq 1$, by the
definition of $u_a$ in (2.7), we have that
$u_a(t_\ve)=\ds\f{\ve}{r(t_\ve)}V(t_\ve)=\ds\f{\ve}{r(t_\ve)}
V\big(r(t_\ve) -t_\ve,\mu\tau_0\big)$ holds. In addition, it follows
from Proposition 2.4 that
$$
|Z^\al(u-u_a)|\leq C_{\alpha, \mu}\varepsilon^{3/2}\,|\ln\ve|
\left(1+t\right)^{-1}\left(1+|t-r|\right)^{1/2}.
$$
Therefore,
$$
\align
w_1(t_\ve)
&=(\p_tu+r\p_{tr}u)(t_\ve)-\left[c(2\p_ru+r\p_r^2u)\right](t_\ve)\\
&=(\p_tu_a+r\p_{tr}u_a)(t_\ve)-\left[c(2\p_ru_a+r\p_r^2u_a)\right](t_\ve)+o(\ve)\\
&=-2\ve\p_q^2V(t_\ve)+\f{\ve^2}{1+t}\p_{q\tau}V(t_\ve)+o(\ve)\\
&=-2\ve\p_q^2V(t_\ve)+o(\ve)
\endalign
$$
which yields $\widehat{w}(t_\ve)=2\ve\p_q^2V(t_\ve)+o(\ve)$.

By Lemma 3.3, we obtain
$$
\biggl(\int_{t_\ve}^{T_\ve}a_0(t)\,dt\biggr)
\exp\biggl(-\int_{t_\ve}^{T_\ve}|a_1(t)|\,dt\bigr)<
\left(\widehat{w}(t_\ve)-K\right)^{-1},
$$
and then
$$
\biggl(\int_{t_\ve+\ve/2}^{T_\ve}\f{1}{4r(t)c(t)}\exp
\biggl(\ds\int_{t_\ve}^{t_\ve+\ve/2}\f{L_2u}{4c^2}(s,r(s))ds\biggr)\,dt\biggr)
\exp\biggl(-\int_{t_\ve}^{T_\ve}|a_1(t)|\,dt\biggr)
<\left(\widehat{w}(t_\ve)-K\right)^{-1}.
$$
It follows
$$
\big(\ln (T_{\ve}+1)-\ln(e^{\mu\tau_0/\ve}+\ve/2)\big)(1+O(\ve))
<4(2\ve\p_q^2V(t_\ve)+o(\ve))^{-1}\exp(E\ve^2).
$$

In order to obtain (3.27), we need to estimate $V(t_\ve)$. Set
$\bar{q}(\tau)=r(t)-t$. Then
$$\cases
\ds\f{d\bar{q}}{d\tau}=\f{t+1}{r(t)(1+c(t))}(V-\p_qV)(\bar{q},\tau)
+O(\ve^{1/2}|\ln\ve|),\\
\bar{q}(0)=\rho_0.
\endcases
$$
Recalling that, in Lemma 2.1, $q(\tau,\rho_0)$ is the characteristics
of the approximate equation which satisfies $q(0,\rho_0)=\rho_0$, and
then comparing $\bar{q}$ with $q$ yields
$$\cases
\ds\f{d(\bar{q}-q)}{d\tau}=\biggl(\f{2(t+1)}{r(t)(1+c(t))}-1\biggr)
\f{V-\p_qV}{2}(\bar{q},\tau)
+\f{(V-\p_qV)(\bar{q},\tau)}{2}-\f{(V-\p_qV)(q,\tau)}{2}+O\left(\ve^{1/2}|\ln\ve|\right)\\
\q\q=\ds\f{(V-\p_qV)(\bar{q},\tau)}{2}-\f{(V-\p_qV)(q,\tau)}{2}
+O\left(\ve^{1/2}|\ln\ve|\right),\\
\bar{q}(0)-q(0,\rho_0)=0.
\endcases
$$
By Gronwall's inequality, we have $\bar{q}(\tau)-q(\tau,\rho_0)
=O\left(\ve^{1/2}|\ln\ve|\right)$ when $\tau\leq R$. Thus,
$\p_q^2V(t_\ve)=\p_q^2V(q(\mu\tau_0,\rho_0),\rho_0)
+O\left(\ve^{1/2}|\ln\ve|\right)$. From the expression for $V$ in
Lemma 2.1, we have that
$\p_q^2V(q(\tau,s),\tau)=\ds\f{F''_0(s)}{\p_sq(\tau,s)}$, and consequently
$\p_q^2V(t_\ve)=\ds\f{F''_0(\rho_0)}{\p_sq(\mu\tau_0,\rho_0)}
+O\left(\ve^{1/2}|\ln\ve|\right)$.  Moreover, one has
$F''_0(\rho_0)>0$ by Remark~1.2. Therefore,
$$
\overline{\lim_{\ve\rightarrow0}}\ve\ln (T_{\ve}+1)-\mu\tau_0\leq
\f{2\p_sq(\mu\tau_0,\rho_0)}{F''_0(\rho_0)},
$$
and as $\mu\rightarrow 1$,
$$
\overline{\lim_{\ve\rightarrow0}}\ve\ln T_{\ve}\leq \tau_0+\f{2\p_sq(\tau_0,\rho_0)}{F''_0(\rho_0)}=\tau_0.
$$
This finishes the proof of (3.27). \qed
\enddemo

\demo{Proof of Theorem 1.1}
Follows now directly from (2.31) and (3.27). \qed
\enddemo



\enddocument